%% file: lame1.tex
\documentclass[12pt]{article}
\usepackage{epsfig}
\input cyracc.def
\font\tc=wncyr10

\def\R{{\mathbf{R}}}
\title{On the hyperbolic metric of the complement of a
rectangular lattice}
\author{Alexandre Eremenko\thanks{Supported by NSF grant DMS-1067886.}}
\date{October 2011}
\begin{document}
\maketitle
\begin{abstract} The density of the hyperbolic metric on the complement
of a rectangular lattice is investigated. The question is related to
conformal mapping of symmetric circular quadrilaterals with all zero angles.

MSC: 30C20.

Keywords: hyperbolic metric, conformal mapping, Lam\'e equation, 
accessory parameter, Landau's constant.
\end{abstract}

\noindent
{\bf 0. Introduction.} A famous theorem of Landau says that there exists
$L>0$ with the property that for all functions
$f$ analytic in the unit disc, such that $|f'(0)|=1$, the image
contains a disc of radius $L-\epsilon$, for every $\epsilon>0$.

The largest constant $L$ for which this is true is called the Landau constant,
and its exact value is unknown. It is natural to conjecture that the extremal
functions are universal coverings of complements to lattices in the plane.

For lattices with the property that the largest disc in the complement
has radius $1$, one can consider the problem of maximizing
$|f'(0)|$ over all analytic functions mapping the unit disc to the complement
of the lattice. Albert Baernstein II and J. Vinson \cite{Bae}
found that the hexagonal lattice and
the universal covering with the property that $f(0)$ is a center of
the complement give a strict local maximum in this problem.

This paper originated from an attempt of the author to prove a similar
result for the class of rectangular lattices. A natural conjecture is 
that the square lattice must be extremal, and that for the extremal function
$f(0)$ is a center of the complement. This is still unproved.
However, the author feels that even the computation of the hyperbolic
metric at the center of the complement of a rectangular lattice
is an interesting
problem.  We will see that it is related to conformal mapping of certain
symmetric circular quadrilaterals, whose all angles are zero, the problem
which attracted attention of mathematicians and physicists
during the whole XX century \cite{Fok,Hilbert,Keen, Klein,Nehari,Venkov}.
\vspace{.1in}

We consider a rectangular lattice $\{ 2m\omega+2n\omega'\}$ where
$\omega=(\ln K)/2,$
$K>1$, and $\omega'=\pi i$. Let $f$ be a universal cover of the complement
of this lattice by the unit disc such that $f(0)=\omega+\omega'$
is a center. 
We are interested in the quantity $A(\ln K)=|f'(0)|$
as a function of $\ln K$.
\vspace{.1in}

In our use of the standard notation
of the theory of elliptic functions we follow \cite{HC} (see also
\cite{Akh}):  
$\tau=\omega'/\omega,\quad h=e^{\pi i\tau};
$
$\theta_j(\zeta)$ is the $j$-th theta-function, 
$\theta_j=\theta_j(0)$, and Jacobi's Modular Function is $\kappa^2$.
\vspace{.1in}

\noindent
1. Denote $$k=\frac{\ln K}{2\pi},\quad a(k)=A(2\pi k)=A(\ln K)=|f'(0)|.$$
We may assume that $f$ maps a circular quadrilateral $Q$
(having zero angles, inscribed in the unit circle, symmetric with respect
to reflections in coordinate axes) onto a fundamental rectangle $R$ of the
lattice, such that $f(0)$ is the center of $R$.
Then a simple symmetry and rescaling argument gives the
functional equation
$$a(k^{-1})=k^{-1}a(k), \quad k>0.$$
It is easy to see that in the limit when $k\to 0$, $f$ maps the unit disc
onto the strip $0<|\Im u|<2\pi$, and we obtain
$$a(0)=4.$$
Together with the functional equation this implies
$$a(k)\sim 4k,\quad k\to\infty.$$
Differentiating the functional equation we obtain
$$a'(1)=\frac{1}{2}a(1).$$
In the next two sections we find $a(1)$ and $a'(0)$. In section 4, we derive an
``explicit expression'' for $a(k)$ for all $k$, and in section 5
discuss the computational aspects.

Finding the density of hyperbolic metric
in the complement of a rectangular lattice is equivalent to
finding a conformal map from a rectangle onto a hyperbolic quadrilateral
with zero angles.
This is a classical problem investigated by Hilbert \cite{Hilbert}
and Klein \cite{Klein}, and in modern times in \cite{Keen}.
The mapping satisfies a Schwarz
differential equation related to a Lam\'e equation.
The problem of finding this mapping for a given circular quadrilateral
(not necessarily inscribed in a circle) requires determination of the
so-called accessory parameter which is a solution of a transcendental
equation involving  Hill's determinants. See
\cite{Fok,Nehari,Venkov} for results in this direction.
The only paper I know where the accessory parameter was actually
computed is \cite{Keen} but it does not contain a rigorous
justification of the algorithm.\footnote{The authors say on p. 217:
``It should be emphasized that our remarks about the implicit equations
are purely heuristic and that the actual computation proceeded, as it were,
fortuitously without any a priori justification.}
So we describe in sections 4 and 5 a convergent algorithm for the special
case considered here.
\vspace{.1in}

\noindent
{\bf 2. Finding $a(1)$.} We use the upper half-plane rather than the unit disc.
Let $T$ be the triangle in the upper half-plane with zero angles and
vertices $0,1,\infty$. Quadrilateral $Q_1$ is the union of
$T$ with its reflection about its left vertical side. Now
the ``center'' is at the point $\tau=i$, and we are looking for
$|f_1^\prime(i)|$, where 
$$f_1(\tau)=f(\zeta),$$
$\tau(\zeta)$ is a map prom the unit disc onto the upper halfplane,
$\tau(0)=i$. 
Our function $f_1$ maps $Q_1$ onto a square $R_1$ of side $2\pi$.
We split it into a composition of two functions.
$$f_1=g\circ\kappa^2,$$
where $\kappa^2$ is the Modular Function of Jacobi.
(In \cite{Akh} this function has double notation, sometimes $\lambda$,
sometimes $\kappa^2$).
$\kappa^2$ maps $T$ onto the upper half-plane sending
$(\infty,0,1)$ to $(0,1,\infty)$. Thus $\kappa^2(i)=1/2$.

Our second component is a Schwarz-Christoffel map
$$g(w)=C\int_0^w\zeta^{-3/4}(1-\zeta)^{-3/4}d\zeta,$$
with
$$\displaystyle
C=\frac{2\pi\sqrt{2}}{B(1/4,1/4)}\approx 1.981\, ,$$
where $B$ is Euler's Beta-function.
This $g$ maps the upper half-plane onto the isosceles right triangle with
legs $2\pi$, and $g(1/2)$ is the middle of the hypotenuse. We have
$$g'(1/2)=C\cdot 2^{3/2}=3.887\, .$$

Using the notation from \cite{Akh,HC}, we have
$$\kappa^2=(\theta_2/\theta_3)^4,\quad \mbox{\cite[II,4, \S5]{HC}},$$
\begin{eqnarray*}
\theta_2&=&2h^{1/4}(1+h^2+h^6+h^{12}+\ldots),\\ \\
\theta_3&=&1+2h+2h^4+2h^9+\ldots,\qquad\mbox{(see \cite[II,2, \S6]{HC})},
\end{eqnarray*}
where
$$h=\exp(\pi i\tau).$$
Thus
$$\frac{d}{d\tau}\kappa^2=4\theta_2^3\theta_3^{-5}(\theta_2^\prime\theta_3-
\theta_2\theta_3^\prime)\frac{dh}{d\tau}.$$
and $d\tau/d\zeta=2$ at $z=0$.

For $\zeta=0,\,\tau=i,\, h=e^{-\pi}$, Matlab gives 
$$|f'(0)|=2|(\kappa^2)'(i)|g'(1/2)\approx 7.416.$$
\vspace{.1in}

\noindent
{\bf 3. Finding $a'(0)$.}
Let us make a preliminary map of the unit disc onto the horizontal
strip $|\Im w|<2\pi$, so that our quadrilateral $Q$ is mapped onto the
quadrilateral $Q_2$ bounded by two vertical segments of length $2\pi$
each, and two half-circles, orthogonal to the boundary of the strip,
and zero corresponds to the center of $Q_2$, which is also at $w=0$.
We have
$$\left|\frac{dw}{d\zeta}\right|=4\quad\mbox{at the point}\quad \zeta=0.$$
Suppose that $Q_2$ has width $\epsilon>0$ which is small.
Then $f_2$ maps $Q_2$ into a rectangle $R_2$, bounded by two
vertical segments of length $2\pi$ and two horizontal
segments, and also having small width. It is enough to estimate the width
of this $R_2$ because
$$f_2^\prime(0)\approx(\mbox{width of}\,R_2)/\epsilon.$$
Let us rescale both $Q_2$ and $R_2$ to the width $1$, denoting the rescaled
quadrilateral by $Q_2^\prime$ and $R_2^\prime$, and let $G$ be the conformal
map between them.
Now we are interested in the ratio of the lengths of these quadrilaterals,
which are large.
If we cut $Q_2^\prime$ and $R_2^\prime$ by horizontal segments in the middle,
the lower half of $Q_2^\prime$ will be conformally equivalent to
the lower half of $R_2^\prime$ (as curvilinear quadrilaterals).
Let us compare the restriction of $G$ to these halfs
with the conformal map $F$
of the triangle $T$ (see section 1) onto a vertical halfstrip with
vertices $0,1,\infty$ (and right angles at $0$ and $1$). Such map with
the vertex correspondence $(0,1,\infty)\to(0,1,\infty)$ is given
by 
$$F(\tau)=\frac{1}{\pi}\arccos\frac{2-\kappa^2(\tau)}{\kappa^2(\tau)}.$$
Using the explicit expression for $\kappa^2$ above and putting $\tau=iy$, we obtain
$$F(iy)/i=y-\frac{1}{\pi}\ln 4+o(1),\quad y\to+\infty.$$
This means that $R_2^\prime$ is shorter than $Q_2^\prime$ by 
$$\frac{2}{\pi}\ln 4+o(1).$$
(That $G(iy)=F(iy)+o(1), y\to\infty$ is quite evident, and this is
easy to justify by using extremal length or some other argument).
Passing back to the original $R_2$ and $Q_2$ we obtain
that the width of $R_2$ is
$$\epsilon+\frac{\epsilon^2}{\pi^2}\ln 4+o(\epsilon^2).$$
This means that $a=|f'(0)|=4|f_2^\prime(0)|=
4+4\epsilon\pi^{-2}\ln 4+o(\epsilon),$ that is
$$a'(0)=\frac{4}{\pi^2}\ln 4\approx0.5618.$$
\vspace{.1in}

\noindent
{\bf 4. Representation of $|f'(0)|=a(k)$}
for arbitrary~$k$.

Let $\wp$ be the Weierstrass function of our
rectangular lattice (with real period $2\omega$ and pure imaginary
period $2\omega'$). We set
$$P(\zeta)=\frac{1}{4}(\wp(\zeta+\omega+\omega')-e_2).$$
Then $P$ is real on both real and imaginary axis, in fact it maps
the rectangle $R_0=\{ 0<\Re u<\omega, \,0<\Im u<\omega'/i\}$ (one quarter
of the fundamental rectangle) onto
the lower half-plane (see Fig. 1)
$P$ is holomorphic in the closure of $R_0$, except at
one point $\omega+\omega'$ where it has a pole of second order.
\vspace{.2in}

\begin{center}
\input{lam1.tex}
\vspace{.2in}

Figure 1. Conformal mapping by $P$.
\end{center}
\vspace{.2in}
Consider the differential equation
\def\ds{\displaystyle}
\begin{equation}
\label{system}
\frac{d^2 w}{d\zeta^2}+P(\zeta)w=\lambda w,\end{equation}
where $\lambda=\lambda(\omega,\omega')$ is a real parameter to be specified
later.

Let $\lambda^-$ be the largest eigenvalue of (\ref{system})
with boundary conditions  $w'(0)=w'(\omega)=0$. Taking into account that
$P(\zeta)<0$ for $\zeta\in[0,\omega]$ (see Fig. 1), we conclude that
\begin{equation}
\label{sturm2}
\lambda^-\leq 0.
\end{equation}
Putting $\zeta=it$ in (\ref{system}) we obtain
\begin{equation}
\label{rotated}
\frac{d^2w}{dt^2}=(P(it)-\lambda)w.
\end{equation}
Let $\lambda^+$ be the smallest eigenvalue of (\ref{rotated}) with
the boundary conditions $w'(0)=w'(\omega'/i)=0$. Using
$P(\zeta)>0$ for $\zeta\in (0,\omega']$ (see Fig. 1),
we obtain
\begin{equation}
\label{sturm1}
\lambda^+\geq 0.
\end{equation}
Let $c$ and $s$ be two solutions of (\ref{system}) defined by
$$\left(\begin{array}{cc} c(0)&s(0)\\
                          c'(0)&s'(0)\end{array}\right)=
\left(\begin{array}{cc} 1&0\\
                          0&1\end{array}\right).$$
Then $c(\sigma),s(\sigma)$ are real for $\sigma\in[0,\omega]$,
while $c(it)$ is real and $s(it)$ is pure imaginary for $t\in[0,\omega']$.
We introduce four real quantities:
\begin{eqnarray}\label{quantities}
m&=&(s(\omega)/c(\omega)+s'(\omega)/c'(\omega))/2,\\
&&\nonumber\\
r&=&(s(\omega)/c(\omega)-s'(\omega)/c'(\omega))/2,\\
&&\nonumber\\
m_1&=&-i(s(\omega')/c(\omega')+s'(\omega')/c'(\omega'))/2,\\
&&\nonumber\\
r_1&=&-i(s(\omega')/c(\omega')-s'(\omega')/c'(\omega'))/2.
\end{eqnarray}

\noindent
{\bf Proposition.} 
{\em There is a unique $\lambda\in[\lambda^-,\lambda^+]$ such that
the condition
\begin{equation}
\label{cond1}
\frac{m_1^2}{\sqrt{m_1^2+m^2}}=r_1
\end{equation}
is satisfied. Then we have}
\begin{equation}
\label{cond2}
a(k)=\frac{mm_1}{\sqrt{m^2+m_1^2}}.
\end{equation}

{\em Proof of the Proposition}. Consider the Schwarz differential equation
associated with (\ref{system}):
\begin{equation}
\label{schwarz}
\frac{F'''}{F'}-\frac{3}{2}\left(\frac{F''}{F'}\right)^2=2(P-\lambda).
\end{equation}
It is well known (and easy to verify) that every solution of
(\ref{schwarz}) is a ratio of two linearly independent solutions of
(\ref{system}). So we put $F=s/c$.
\vspace{.1in}

\noindent
4.1.
Then $F$ is locally univalent, real on the real axis, and pure imaginary
on the imaginary axis.
\vspace{.1in}

\noindent
4.2. 
We claim that for $\lambda\in [\lambda^- ,\lambda^+],$
\begin{equation}
c(\zeta)\neq 0\quad\mbox{when}\quad \zeta\in [0,\omega]\cup [0,\omega'].
\end{equation}
This follows from the Sturm Comparison Theorem.

Thus $F$ is holomorphic on $[0,\omega]\cup[0,\omega']$.
\vspace{.1in}

\noindent
4.3.
By the Symmetry principle, 
\begin{equation}
\label{symmetry}
F(\overline{\zeta})=\overline{F(\zeta)}\quad\mbox{and}\quad
F(-\overline{\zeta})=-\overline{F(\zeta)}.
\end{equation}
It is well-known that $F$ maps the right vertical side $L'$
and the top horizontal side $L$ of the rectangle
$R_0$ onto some arcs of circles on the Riemann sphere.
We recall a simple proof of this.
Consider, for example the right vertical side $L'$.
As the right hand side of (\ref{schwarz}) has period $2\omega$,
and all solutions of (\ref{schwarz}) are fractional linear transformations
of each other, we conclude that $F(\zeta-2\omega)=\lambda\circ F(\zeta)$.
Using  $L'-2\omega=-\overline{L'}$ and (\ref{symmetry}),
we obtain
$$-\overline{F(L')}=\lambda\circ F(L').$$
Thus $F(L')$ is fixed by
an anticonformal involution, so it is an arc of a circle.
Similar argument applies to $F(L)$. 
\vspace{.1in}

\noindent
4.4.
It is known that the arcs $F(L)$ and $F(L')$ have exactly one
common
point $z_0$ and they are tangent at this point. We recall how this is proved.
The coefficient $P$ in (\ref{system}) 
has one pole in the closure of $R_0$, namely at the point $\omega+\omega'$,
and the Laurent series at this pole has the form 
$$P(\omega+\omega'+\zeta)=\frac{1}{4\zeta^2}+\ldots .$$
It follows that (\ref{system}) has two linearly independent solutions
near $\omega+\omega'$ of the form
$$\phi(\omega+\omega'+\zeta)=\zeta^{1/2}+\ldots\quad\mbox{and}\quad
\psi(\omega+\omega'+\zeta)=\zeta^{1/2}\ln\zeta+\ldots$$
(see, for example, \cite{I}).
Thus $F$ (which is a fractional-linear transformation of $\phi/\psi$)
has a limit $z_0$ as $\zeta\to \omega+\omega'$ from inside of $R_0$, and
the angle between $F(L)$ and $F(L')$ at $z_0$ is zero. As
$F(L)$ and $F(L')$ are arcs of circles, $z_0$ is the only point of
intersection of these circles.
\vspace{.1in}

\noindent
4.5.
We have shown in 4.1-4.4 that $F$ maps the boundary $\partial R_0$ locally
univalently onto a Jordan curve (consisting of a segment of the real
axis, a segment of the imaginary axis, and two arcs of circles
tangent to each other, one perpendicular to the real axis
another to the imaginary axis).
\vspace{.1in}

\begin{center}
\input{lam22.tex}
\vspace{.2in}
Figure 2. Conformal mapping by $F$.
\end{center}
\vspace{.2in}

We claim that right circle has center at the point $m$ and radius $r$,
while the top circle has center $im_1$ and radius $r_1$ (the definitions
of these four numbers are given in (\ref{quantities}).

Let us verify the claim for the right circle.
It is clear that the matrix
$$\left(\begin{array}{cc}s(\omega)&s'(\omega)\\
           c(\omega)&c'(\omega)\end{array}\right)$$
is real. Introducing the pair of solutions
$(\phi,\psi)$ of the equation (\ref{rotated}), normalized by
$$\left(\begin{array}{cc}\phi(\omega)&(d\phi/dt)(\omega)\\
                         \psi(\omega)&(d\psi/dt)(\omega)\end{array}\right)=
\left(\begin{array}{cc}1&0\\0&1\end{array}\right),\quad\mbox{where}\quad
\zeta=it,$$
we see that $\phi$ and $\psi$ are real on the vertical line
$\{\omega+it:t\in\R\}$ (as solutions of real differential
equation (\ref{rotated}) with real initial conditions),
and conclude that
\begin{eqnarray*}
s&=&s(\omega)\phi+is'(\omega)\psi,\\
c&=&c(\omega)\phi+ic'(\omega)\psi,
\end{eqnarray*}
and thus the image $F(L')$ belongs to the circle
$$
\left\{\frac{s(\omega)+s'(\omega)ix}{c(\omega)+c'(\omega)ix}:x\in\R\right\},$$
whose center is at the point $m$ and radius is $r$.

The verification for the other circle is similar.

Now we choose our parameter $\lambda$ in such a way that the
common tangent line to the two circles passes through the origin.
Elementary geometry shows that this happens iff the condition (\ref{cond1})
is satisfied, see Figure 3.
\begin{center}
\input{lam3.tex}
\vspace{.2in}

Figure 3. Choosing the accessory parameter.
\end{center}
\vspace{.2in}

When $\lambda\to\lambda^-$, the right circle becomes a vertical line,
and when $\lambda\to\lambda^+$ the top circle becomes a horizontal line.
It follows that the quantity $e=m_1^2/\sqrt{m^2+m_1^2}-r_1$ changes sign
when $\lambda$ changes on the interval $[\lambda^-,\lambda^+]$.
It follows that there exists $\lambda$ such that condition (\ref{cond1})
is satisfied.
In fact this value of $\lambda$ is unique but we do not use this additional
information.

For this value of $\lambda$, $F$ maps our rectangle $R_0$ onto a quadrilateral
$Q_0$ bounded by two straight segments $F([0,\omega])\subset\R$,
$F([0,\omega'])\subset i\R$ and two arcs of circle tangent at the point $z_0$.
Reflecting this
quadrilateral $Q_0$ three times with respect to the coordinate axes
we obtain
a circular quadrilateral symmetric with respect to the coordinate axes,
inscribed in the circle $\{ z:|z|=|z_0|\}$ and sides perpendicular
to this circle. Thus $F$ is the inverse to the universal
covering of the complement of our lattice by the disc 
$\{ z:|z|<|z_0|\}$, $F(0)=0$, $F'(0)=1$,
and it remains to verify that the radius $|z_0|$ of
this disc is given by the formula (\ref{cond2}). This is clear from Figure 3.
\vspace{.1in}

\noindent
{\bf 5. Remarks on computation}. 
Our method of computation is based on the previous section.
We are solving the Lam\'e equation (\ref{system}). We represent $P$ in
terms of theta-functions or in terms of elliptic functions of Jacobi 
(see \cite{HC,Akh}):
$$
P(\zeta)=-\frac{{\theta_1^\prime}^2\,\theta_1^2(\zeta/(2\omega))}{
                  16\omega^2\,\theta_3^2\,\theta_3^2(\zeta/(2\omega))}\\
$$
where either product or series representations can be used for theta
functions, they converge very well.
Numerical solving of the equations (\ref{system}) and (\ref{rotated})
on the real line involves only
computations with real numbers.

As Matlab does not have standard routines for theta functions,
we can use an expression of $\wp$ in terms of Jacobi elliptic
functions. Matlab uses AGM to compute them which is probably as effective
as theta functions.

To find the accessory parameter $\lambda$ we use a ``shooting method'',
dissecting the interval $[\lambda^-,\lambda^+]$ dyadically.

The result is that $a(1)\approx 1.4163$ and the graph
of $|f'(0)|$ for the rectangular lattice
with sides $2\omega,2\omega',\;4\omega^2+4{\omega^\prime}^2=1$,
as a function of $2\omega$ is this:
\begin{center}
\epsfxsize=5.2in
\centerline{\epsffile{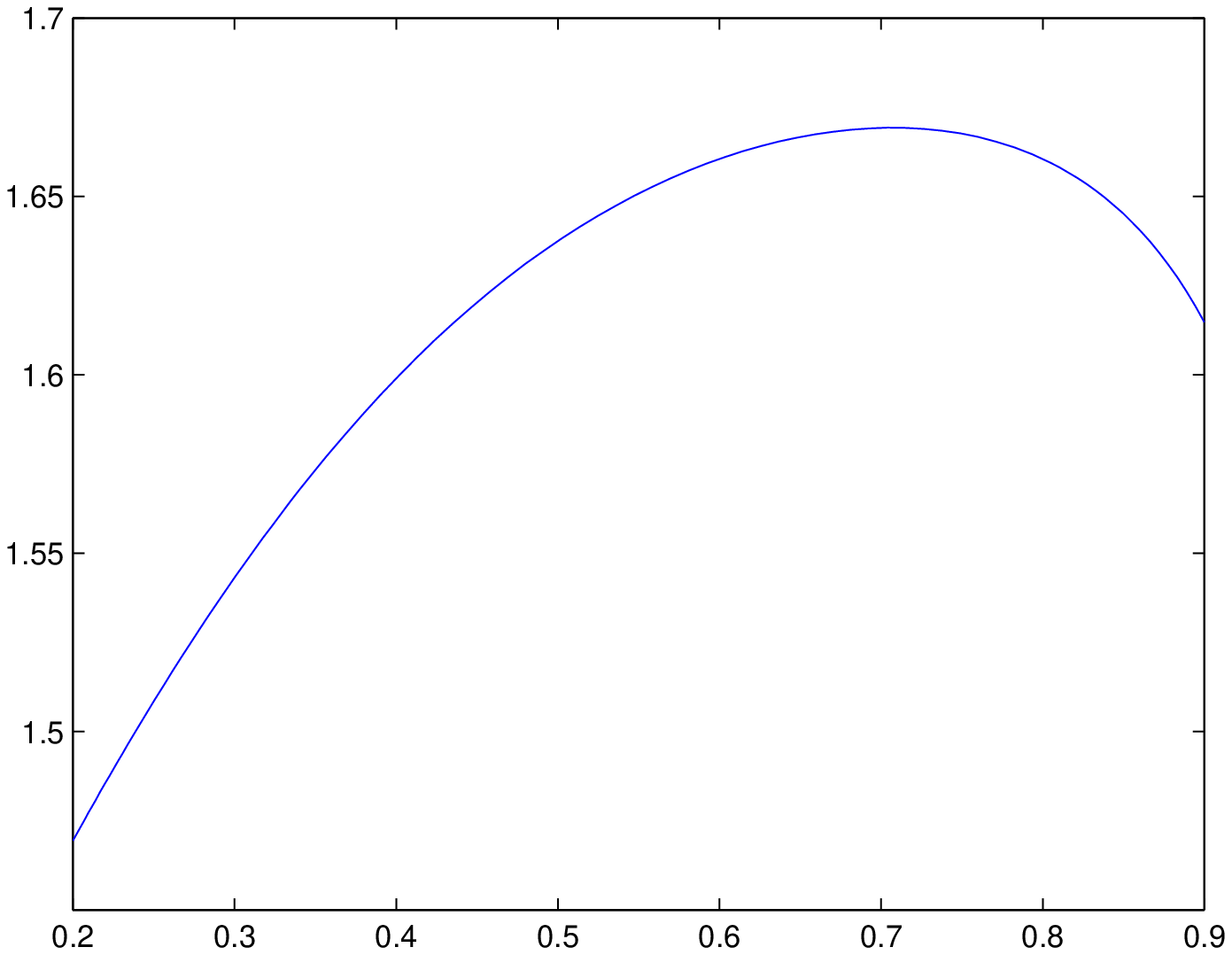}}
\end{center}
We see a maximum at the point $1/\sqrt{2}$ and this maximum
equals $1.6693$.

{\em Purdue University, West Lafayette IN 47907

eremenko@math.purdue.edu}
\end{document}

%% file: lam1.tex
\begin{picture}(0,0)%
\includegraphics{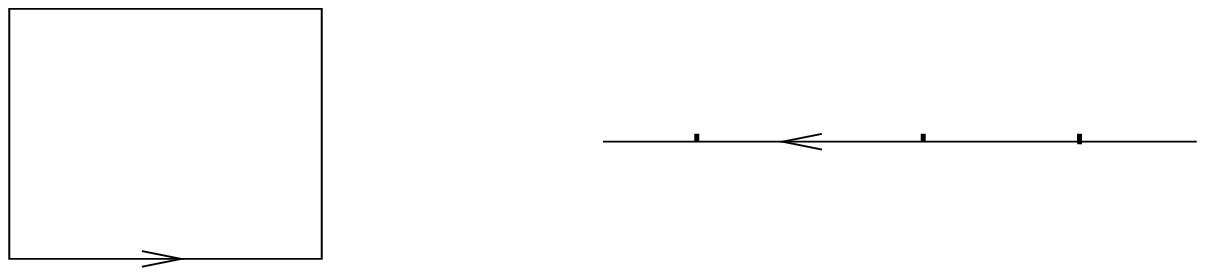}%
\end{picture}%
\setlength{\unitlength}{1973sp}%
\begingroup\makeatletter\ifx\SetFigFont\undefined%
\gdef\SetFigFont#1#2#3#4#5{%
  \reset@font\fontsize{#1}{#2pt}%
  \fontfamily{#3}\fontseries{#4}\fontshape{#5}%
  \selectfont}%
\fi\endgroup%
\begin{picture}(11722,4413)(1201,-5986)
\put(4501,-4936){\makebox(0,0)[lb]{\smash{\SetFigFont{12}{14.4}{\rmdefault}{\mddefault}{\updefault}
$\omega$
}}}
\put(1201,-4861){\makebox(0,0)[lb]{\smash{\SetFigFont{12}{14.4}{\rmdefault}{\mddefault}{\updefault}
$0$
}}}
\put(1201,-1936){\makebox(0,0)[lb]{\smash{\SetFigFont{12}{14.4}{\rmdefault}{\mddefault}{\updefault}
$\omega'$
}}}
\put(4051,-1861){\makebox(0,0)[lb]{\smash{\SetFigFont{12}{14.4}{\rmdefault}{\mddefault}{\updefault}
$\omega+\omega'$
}}}
\put(4726,-2311){\makebox(0,0)[lb]{\smash{\SetFigFont{12}{14.4}{\rmdefault}{\mddefault}{\updefault}
$(\infty)$
}}}
\put(2626,-3361){\makebox(0,0)[lb]{\smash{\SetFigFont{12}{14.4}{\rmdefault}{\mddefault}{\updefault}
$R_0$
}}}
\put(10051,-3811){\makebox(0,0)[lb]{\smash{\SetFigFont{12}{14.4}{\rmdefault}{\mddefault}{\updefault}
$(0)$
}}}
\put(2401,-5911){\makebox(0,0)[lb]{\smash{\SetFigFont{12}{14.4}{\rmdefault}{\mddefault}{\updefault}
$\zeta$
}}}
\put(9226,-5986){\makebox(0,0)[lb]{\smash{\SetFigFont{12}{14.4}{\rmdefault}{\mddefault}{\updefault}
$P(\zeta)$
}}}
\put(11326,-3061){\makebox(0,0)[lb]{\smash{\SetFigFont{12}{14.4}{\rmdefault}{\mddefault}{\updefault}
$\frac{e_1-e_2}{4}$
}}}
\put(7651,-3061){\makebox(0,0)[lb]{\smash{\SetFigFont{12}{14.4}{\rmdefault}{\mddefault}{\updefault}
$\frac{e_3-e_2}{4}$
}}}
\put(10276,-3061){\makebox(0,0)[lb]{\smash{\SetFigFont{12}{14.4}{\rmdefault}{\mddefault}{\updefault}
$0$
}}}
\put(7726,-3811){\makebox(0,0)[lb]{\smash{\SetFigFont{12}{14.4}{\rmdefault}{\mddefault}{\updefault}
$(\omega)$
}}}
\put(11401,-3811){\makebox(0,0)[lb]{\smash{\SetFigFont{12}{14.4}{\rmdefault}{\mddefault}{\updefault}
$(\omega')$
}}}
\end{picture}

%% file: lam22.tex
\begin{picture}(0,0)%
\includegraphics{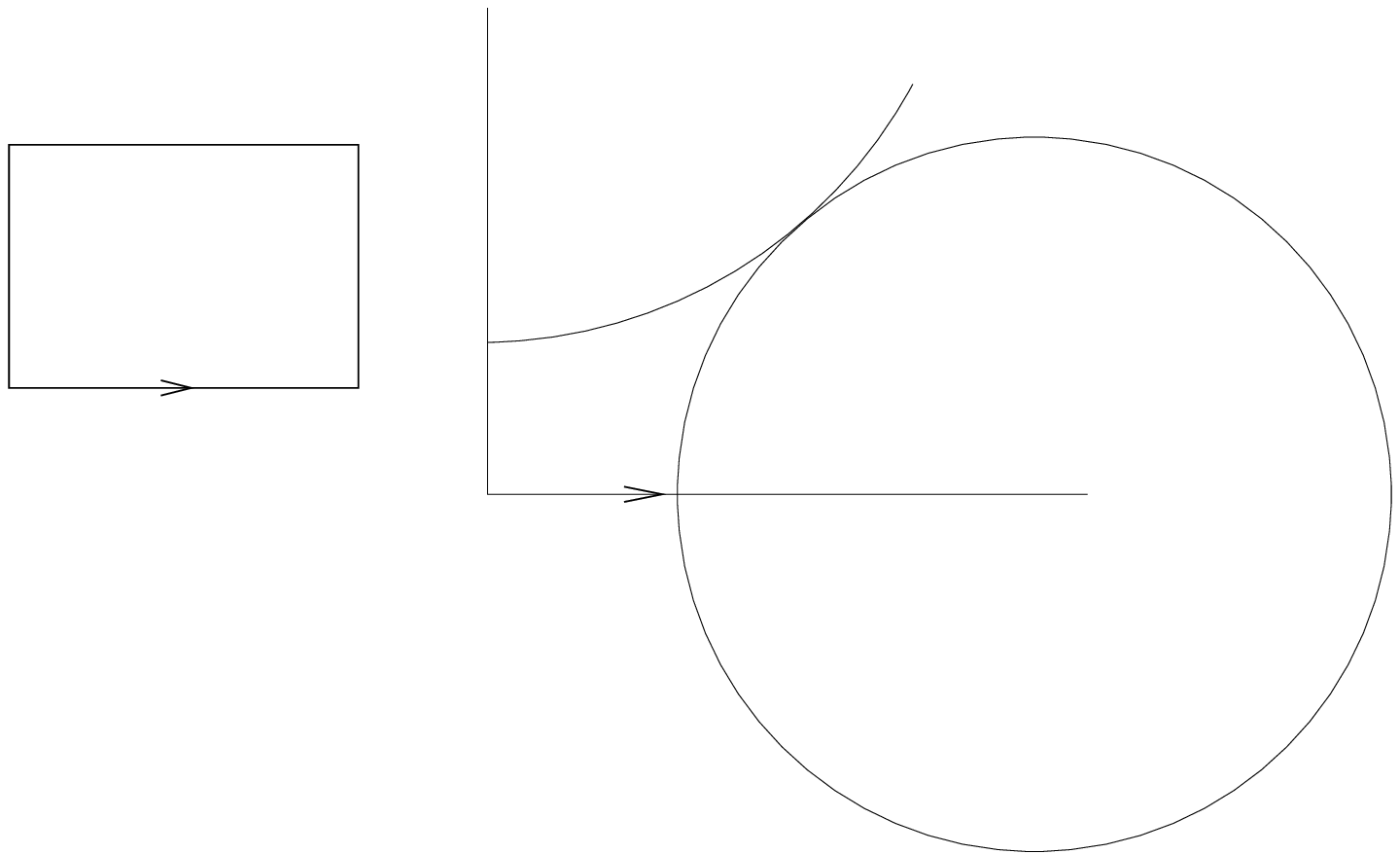}%
\end{picture}%
\setlength{\unitlength}{1973sp}%
\begingroup\makeatletter\ifx\SetFigFont\undefined%
\gdef\SetFigFont#1#2#3#4#5{%
  \reset@font\fontsize{#1}{#2pt}%
  \fontfamily{#3}\fontseries{#4}\fontshape{#5}%
  \selectfont}%
\fi\endgroup%
\begin{picture}(13958,8345)(301,-9144)
\put(301,-2086){\makebox(0,0)[lb]{\smash{\SetFigFont{9}{10.8}{\rmdefault}{\mddefault}{\updefault}
$\omega'$
}}}
\put(4126,-3436){\makebox(0,0)[lb]{\smash{\SetFigFont{9}{10.8}{\rmdefault}{\mddefault}{\updefault}
$L'$
}}}
\put(376,-4861){\makebox(0,0)[lb]{\smash{\SetFigFont{9}{10.8}{\rmdefault}{\mddefault}{\updefault}
$0$
}}}
\put(2026,-3436){\makebox(0,0)[lb]{\smash{\SetFigFont{9}{10.8}{\rmdefault}{\mddefault}{\updefault}
$R_0$
}}}
\put(2101,-2011){\makebox(0,0)[lb]{\smash{\SetFigFont{9}{10.8}{\rmdefault}{\mddefault}{\updefault}
$L$
}}}
\put(4051,-4861){\makebox(0,0)[lb]{\smash{\SetFigFont{9}{10.8}{\rmdefault}{\mddefault}{\updefault}
$\omega$
}}}
\put(5101,-5911){\makebox(0,0)[lb]{\smash{\SetFigFont{9}{10.8}{\rmdefault}{\mddefault}{\updefault}
$0$
}}}
\put(10501,-5461){\makebox(0,0)[lb]{\smash{\SetFigFont{9}{10.8}{\rmdefault}{\mddefault}{\updefault}
$x_0$
}}}
\put(8251,-2686){\makebox(0,0)[lb]{\smash{\SetFigFont{9}{10.8}{\rmdefault}{\mddefault}{\updefault}
$z_0$
}}}
\put(7426,-4786){\makebox(0,0)[lb]{\smash{\SetFigFont{9}{10.8}{\rmdefault}{\mddefault}{\updefault}
$F(L')$
}}}
\put(6001,-3811){\makebox(0,0)[lb]{\smash{\SetFigFont{9}{10.8}{\rmdefault}{\mddefault}{\updefault}
$F(L)$
}}}
\put(6451,-7336){\makebox(0,0)[lb]{\smash{\SetFigFont{9}{10.8}{\rmdefault}{\mddefault}{\updefault}
$F(\zeta)$
}}}
\put(1651,-7336){\makebox(0,0)[lb]{\smash{\SetFigFont{9}{10.8}{\rmdefault}{\mddefault}{\updefault}
$\zeta$
}}}
\end{picture}

%% file: lam3.tex
\begin{picture}(0,0)%
\includegraphics{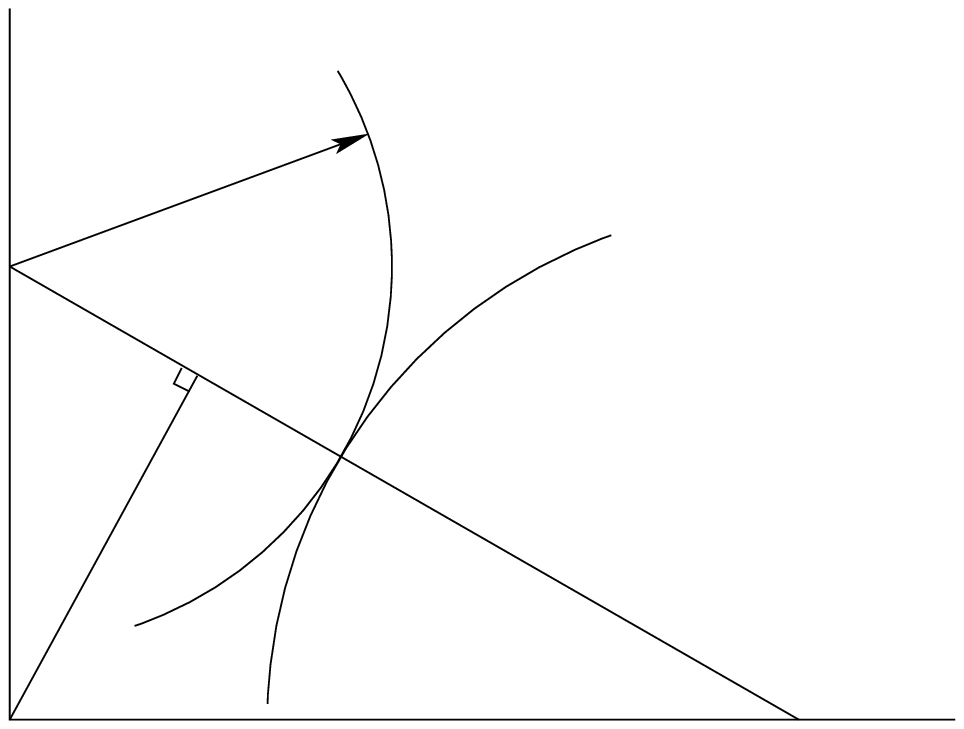}%
\end{picture}%
\setlength{\unitlength}{1973sp}%
\begingroup\makeatletter\ifx\SetFigFont\undefined%
\gdef\SetFigFont#1#2#3#4#5{%
  \reset@font\fontsize{#1}{#2pt}%
  \fontfamily{#3}\fontseries{#4}\fontshape{#5}%
  \selectfont}%
\fi\endgroup%
\begin{picture}(9997,6997)(1201,-7111)
\put(9601,-6736){\makebox(0,0)[lb]{\smash{\SetFigFont{12}{14.4}{\rmdefault}{\mddefault}{\updefault}
$m$
}}}
\put(1201,-2611){\makebox(0,0)[lb]{\smash{\SetFigFont{12}{14.4}{\rmdefault}{\mddefault}{\updefault}
$im_1$
}}}
\put(3376,-1936){\makebox(0,0)[lb]{\smash{\SetFigFont{12}{14.4}{\rmdefault}{\mddefault}{\updefault}
$r_1$
}}}
\put(4501,-3886){\makebox(0,0)[lb]{\smash{\SetFigFont{12}{14.4}{\rmdefault}{\mddefault}{\updefault}
$e$
}}}
\put(1651,-7111){\makebox(0,0)[lb]{\smash{\SetFigFont{12}{14.4}{\rmdefault}{\mddefault}{\updefault}
$0$
}}}
\end{picture}